\documentclass[10pt,twoside]{amsart}
\usepackage{amssymb}
\usepackage{float}
\usepackage{amsfonts}
\usepackage[centertags]{amsmath}
\usepackage{amsthm}
\usepackage{newlfont}
\usepackage{array}

\date{}

\begin{document}

\centerline{\Large{\bf Fixed point results for a new metric space}}

\centerline{}

\centerline{\bf {Murat Kiri\c{s}ci}*, \bf {Necip \c{S}im\c{s}ek} and \bf {Mahmut Akyi\u{g}it} }

\centerline{}

\centerline{a. Department of Mathematical Education, Istanbul University-Cerrahpa\c{s}a, Vefa, 34470, Fatih, Istanbul, Turkey}

\centerline{e-mail: mkirisci@hotmail.com}

\centerline{b. Department of Mathematics, Istanbul Commerce University, Istanbul, Turkey}

\centerline{e-mail: necipsimsek@hotmail.com}

\centerline{}

\centerline{c. Department of Mathematics, Sakarya University, Sakarya, Turkey}

\centerline{e-mail: makyigit@sakarya.edu.tr}

\centerline{}

\newtheorem{theorem}{\quad Theorem}

\newtheorem{definition}[theorem]{\quad Definition}

\newtheorem{corollary}[theorem]{\quad Corollary}

\newtheorem{proposition}[theorem]{\quad Proposition}

\newtheorem{lemma}[theorem]{\quad Lemma}

\newtheorem{example}[theorem]{\quad Example}

\newtheorem*{remark}{Remark}

\centerline{}
{\textbf{Abstract:}
In this paper, we introduce the neutrosophic contractive  and neutrosophic mapping.
We establish some results on fixed points of a neutrosophic mapping.
\\

\centerline{}



\section{Introduction}

Fuzzy Sets (FSs) put forward by Zadeh \cite{Zadeh} has influenced deeply all the scientific fields since the publication of the paper.
It is seen that this concept, which is very important for real-life situations, had not enough solution to some problems in time.
New quests for such problems have been coming up. Atanassov \cite{Atan} initiated Intuitionistic fuzzy sets (IFSs) for such cases.
Neutrosophic set (NS) is a new version of the idea of the classical set which is defined by Smarandache \cite{Smar}.
Examples of other generalizations are FS \cite{Zadeh} interval-valued FS \cite{Turksen}, IFS \cite{Atan},  interval-valued IFS \cite{AtanGarg},
the sets paraconsistent, dialetheist, paradoxist, and tautological \cite{Smar0}, Pythagorean fuzzy sets \cite{Yager} .\\

Using the concepts Probabilistic metric space and fuzzy, fuzzy metric space (FMS) is introduced in \cite{KraMic}.
Kaleva and Seikkala \cite{KalSei} have defined the FMS as a distance between two points to be a non-negative fuzzy number.
In  \cite{GeoVee} some basic properties of FMS studied and the Baire Category Theorem for FMS proved.
Further, some properties such as separability, countability are given and Uniform Limit Theorem is proved in \cite{GeoVee2}.
Afterward, FMS has used in the applied sciences such as fixed point theory, image and signal processing, medical imaging, decision-making et al.
After defined of the intuitionistic fuzzy set (IFS), it was used in all areas where FS theory was studied.
Park \cite{Park} defined IF metric space (IFMS), which is a generalization of FMSs.
Park used George and Veeramani's \cite{GeoVee} idea of applying t-norm and t-conorm to the FMS meanwhile defining IFMS and studying its basic features.\\

Fixed point theorem for fuzzy contraction mappings is initiated by Heilpern \cite{Heilpern}.
Bose and Sahani \cite{Bose} extended the Heilpern’s study. Alaca et al. \cite{Alaca}  are given fixed point theorems related to intuitionistic fuzzy metric spaces(IFMSs).
Fixed point results for fuzzy metric spaces and IFMSs are studied by many researchers \cite{Gregori}, \cite{Imdad}, \cite{Mihet}, \cite{Turk}, \cite{Hussain}.\\

In this paper, fixed point results for NMSs are given.

\section{Preliminaries}

Some definitions related to the fuzziness, intuitionistic fuzziness and neutrosophy are given as follows:\\

The fuzzy subset $F$ of $\mathbb{R}$ is said to be a fuzzy number(FN).
The FN is a mapping $F:\mathbb{R}\rightarrow [0,1]$  that corresponds to each real number $a$ to the degree of membership $F(a)$.

Let $F$ is a FN. Then, it is known that \cite{Kirisci2}

\begin{itemize}
  \item If $F(a_{0})=1$, for $a_{0} \in \mathbb{R}$, $F$ is said to be normal,
  \item If for each $\mu >0$, $F^{-1}\{[0,\tau+\mu)\}$
is open in the usual topology $\forall \tau\in [0,1)$, $F$ is said to be upper semi continuous, ,
  \item  The set $[F]^{\tau}=\{a\in\mathbb{R}:F(a)\geq \tau\}$, $\tau \in [0,1]$ is called $\tau-$cuts of $F$.
\end{itemize}

Choose non-empty set $F$. An IFS in $F$ is an object $U$ defined by
\begin{eqnarray*}
U=\{<a,G_{U}(a),Y_{U}(a)>: a\in F\}
\end{eqnarray*}
where $G_{U}(a):F\rightarrow [0,1]$ and
$Y_{U}(a):F\rightarrow [0,1]$ are functions for all $a\in F$ such that $0\leq G_{U}(a)+Y_{U}(a) \leq 1$ \cite{Atan}.
Let $U$ be an IFN. Then,

\begin{itemize}
  \item an IF subset of the $\mathbb{R}$,
  \item If $G_{U}(a_{0})=1$ and, $Y_{U}(a_{0})=0$ for $a_{0} \in \mathbb{R}$, normal,
  \item  If $G_{U}(\lambda a_{1}+(1-\lambda)a_{2})\geq \min(G_{U}(a_{1}), G_{U}(a_{2}))$, $\forall a_{1},a_{2}\in\mathbb{R}$ and $\lambda\in[0,1]$, then the membership function(MF) $G_{U}(a)$ is called convex,
  \item If $Y_{U}(\lambda a_{1}+(1-\lambda)a_{2})\geq \min(Y_{U}(a_{1}), Y_{U}(a_{2}))$, $\forall a_{1},a_{2}\in\mathbb{R}$ and $\lambda\in[0,1]$, then the nonmembership function(NMF)$Y_{U}(a)$ is concav,
  \item $G_{U}$ is upper semi continuous and $Y_{U}$ is lower semi continuous
  \item $supp U=cl(\{a\in F: Y_{U}(a)<1\})$ is bounded.
\end{itemize}

An IFS $U=\{<a,G_{U}(a),Y_{U}(a)>:a\in F\}$ such that $G_{U}(a)$ and $1-Y_{U}(a)$ are FNs, where $(1-Y_{U})(a)=1-Y_{U}(a)$, and $G_{U}(a)+Y_{U}(a)\leq 1$ is called an IFN.\\

Let's consider that $F$ is a space of points(objects). Denote the $G_{U}(a)$ is a truth-MF,
$B_{U}(a)$ is an indeterminacy-MF and $Y_{U}(a)$ is a falsity-MF, where $U$ is a set in $F$ with $a\in F$. Then, if we take $I=]0^{-},1^{+}[$
\begin{eqnarray*}
&&G_{U}(a): F\rightarrow I,\\
&&B_{U}(a): F\rightarrow I, \\
&&Y_{U}(a): F\rightarrow I,
\end{eqnarray*}

There is no restriction on the sum of $G_{U}(a)$, $B_{U}(a)$ and $Y_{U}(a)$. Therefore,
\begin{eqnarray*}
0^{-}\leq \sup G_{U}(a) + \sup B_{U}(a)+ \sup Y_{U}(a) \leq 3^{+}.
\end{eqnarray*}
The set $U$ which consist of with $G_{U}(a)$, $B_{U}(a)$ and $Y_{U}(a)$ in $F$ is called a neutrosophic sets(NS) and can be denoted by
\begin{eqnarray}\label{NS}
U=\{<a,(G_{U}(a),B_{U}(a), Y_{U}(a))>:a\in F, G_{U}(a), B_{U}(a), Y_{U}(a)\in I \}
\end{eqnarray}

Clearly, NS is an enhancement of $[0,1]$ of IFSs.\\

An NS $U$ is included in another NS $V$, ($U\subseteq V$), if and only if,
\begin{eqnarray*}
&&\inf G_{U}(a) \leq \inf G_{V}(a), \quad  \sup G_{U}(a) \leq \sup G_{V}(a),\\
&&\inf B_{U}(a) \geq \inf B_{V}(a), \quad   \sup B_{U}(a) \geq \sup B_{V}(a),\\
&&\inf Y_{U}(a) \geq \inf Y_{V}(a),\quad   \sup Y_{U}(a) \geq \sup Y_{V}(a).
\end{eqnarray*}

for any $a\in F$. However, NSs are inconvenient to practice in real problems.
To cope with this inconvenient situation, Wang et al \cite{Wang} customized NS's definition and single-valued NSs (SVNSs) suggested.

To cope with this inconvenient situation, Wang et al \cite{Wang} customized NS's definition and single-valued NSs suggested. Ye \cite{Ye}, described the notion of simplified NSs, which may be characterized by three real numbers in the $[0,1]$. At the same time, the simplified NSs' operations may be impractical, in some cases \cite{Ye}.  Hence, the operations and comparison way between SNSs and the aggregation operators for simplified NSs are redefined in \cite{Peng}.\\

According to the Ye \cite{Ye}, a simplification of an NS $U$, in (\ref{NS}), is
\begin{eqnarray*}
U=\left\{<a,(G_{U}(a),B_{U}(a), Y_{U}(a))>:a\in F\right\},
\end{eqnarray*}
which called an simplified NS. Especially, if $F$ has only one element $<G_{U}(a),B_{U}(a), Y_{U}(a)>$ is said to be an simplified NN. Expressly, we may see simplified NSs as a subclass of NSs.\\

An simplified NS $U$ is comprised in another simplified NS $V$ ($U\subseteq V$), iff $G_{U}(a)\leq G_{V}(a)$, $B_{U}(a)\geq B_{V}(a)$
and $Y_{U}(a)\geq Y_{V}(a)$ for any $a\in F$. Then, the following operations are given by Ye\cite{Ye}:

\begin{eqnarray*}
U+V&=&\langle G_{U}(a)+G_{V}(a)-G_{U}(a).G_{V}(a), B_{U}(a)+B_{V}(a)-B_{U}(a).B_{V}(a), Y_{U}(a)+Y_{V}(a)-Y_{U}(a).Y_{V}(a)\rangle,\\
U.V&=&\langle G_{U}(a).G_{V}(a), B_{U}(a).B_{V}(a), Y_{U}(a).Y_{V}(a)\rangle,\\
\alpha. U&=&\langle 1-(1-G_{U}(a))^{\alpha}, 1-(1-B_{U}(a))^{\alpha}, 1-(1-Y_{U}(a))^{\alpha}\rangle \quad \quad for \quad \alpha>0,\\
U^{\alpha}&=&\langle G_{U}^{\alpha}(a), B_{U}^{\alpha}(a), Y_{U}^{\alpha}(a)\rangle  \quad \quad for \quad \alpha>0.
\end{eqnarray*}

Triangular norms (t-norms) (TN) were initiated by Menger \cite{Menger}. In the problem of computing the distance between two elements in space, Menger offered using probability distributions instead of using numbers for distance. TNs are used to generalize with the probability distribution of triangle inequality in metric space conditions. Triangular conorms (t-conorms) (TC) know as dual operations of TNs. TNs and TCs are very significant for fuzzy operations(intersections and unions).

\begin{definition} Give
an operation $\circ : [0,1] \times [0,1] \rightarrow [0,1]$.
If the operation $\circ$ is satisfying the following conditions, then it is called that the operation $\circ$ is \emph{continuous TN}(CTN): For $s,t,u,v\in [0,1]$,
\begin{itemize}
\item [i.] $s \circ 1 = s$
\item [ii.] If $s\leq u$ and $t \leq v$, then $s \circ t \leq u \circ v$,
\item [iii.] $\circ$ is continuous,
\item [iv.] $\circ$ is commutative and associative.
\end{itemize}
\end{definition}

\begin{definition} Give
an operation $\bullet : [0,1] \times [0,1] \rightarrow [0,1]$.
If the operation $\bullet$ is satisfying the following conditions, then it is called that the operation $\bullet$ is \emph{continuous TC}(CTC):
\begin{itemize}
\item [i.] $s \bullet 0 = s$,
\item [ii.] If $s\leq u$ and $t \leq v$, then $s \bullet t \leq u \bullet v$,
\item [iii.] $\bullet $ is continuous,
\item [iv.] $\bullet$ is commutative and associative.
\end{itemize}
\end{definition}

Form above definitions, we note that if we choose $0<\varepsilon_{1}, \varepsilon_{2}<1$ for $\varepsilon_{1}>\varepsilon_{2}$, then there exist $0<\varepsilon_{3}, \varepsilon_{4}<0,1$ such that $\varepsilon_{1}\circ \varepsilon_{3}\geq \varepsilon_{2}$, \quad $\varepsilon_{1} \geq \varepsilon_{4}\bullet \varepsilon_{2}$. Further, if we choose $\varepsilon_{5} \in (0,1)$, then there exist $\varepsilon_{6}, \varepsilon_{7} \in (0,1)$ such that $\varepsilon_{6}\circ \varepsilon_{6}\geq \varepsilon_{5}$ and $\varepsilon_{7}\bullet \varepsilon_{7}\leq \varepsilon_{5}$.

\begin{remark}\label{remark:1}\cite{Park}
Take $\circ$ and $\bullet$ are a CTN and CTC, respectively. For $p,s,t,u,v \in [0,1]$,
\begin{itemize}
\item[a.] If $s>t$, then there are $u,v$ such that $s \circ u \geq t$ and $s \geq t \bullet v$.
\item[b.] There are $p,t$ such that $t \circ t \geq s$ and $s\geq p \bullet p$.
\end{itemize}
\end{remark}

\begin{definition}\cite{KirSim}
Take $F$ be an arbitrary set, $V=\mathcal{N}=\{<a, G(a),B(a),Y(a)> : a\in F\}$ be a NS such that $\mathcal{N}: F \times F \times \mathbb{R}^{+} \rightarrow [0,1]$.
Let $\circ$ and $\bullet$ show the continuous TN and continuous TC, respectively.
The four-tuple $(F, \mathcal{N}, \circ, \bullet)$ is called neutrosophic metric space(NMS) when the following conditions are satisfied.
$\forall a,b,c\in F$,

\begin{itemize}
\item [i.] $0 \leq G(a,b,\lambda) \leq 1$,\quad $0 \leq B(a,b,\lambda) \leq 1$, \quad $0 \leq Y(a,b,\lambda) \leq 1$ \quad $\forall \lambda \in \mathbb{R}^{+}$,
\item [ii.] $G(a,b,\lambda)+B(a,b,\lambda)+Y(a,b,\lambda)\leq 3$,   (for $\lambda \in  \mathbb{R}^{+}$),
\item [iii.] $G(a,b,\lambda)=1$ \quad (for $\lambda >0$)  if and only if $a=b$,
\item[iv.] $G(a,b,\lambda)=G(b,a,\lambda)$ \quad (for $\lambda >0$),
\item[v.] $G(a,b,\lambda)\circ G(b,c,\mu)\leq G(a,c,\lambda+\mu)$ \quad $(\forall \lambda, \mu >0)$,
\item[vi.] $G(a,b,.): [0,\infty)\rightarrow [0,1]$ is continuous,
\item[vii.] $lim_{\lambda\rightarrow \infty}G(a,b,\lambda)=1$ \quad $(\forall\lambda >0)$,
\item [viii.] $B(a,b,\lambda)=0$ \quad (for $\lambda >0$)  if and only if $a=b$,
\item[ix.] $B(a,b,\lambda)=B(b,a,\lambda)$ \quad (for $\lambda >0$),
\item[x.] $B(a,b,\lambda)\bullet B(b,c,\mu)\geq B(a,c,\lambda+\mu)$ \quad $(\forall \lambda, \mu >0)$,
\item[xi.] $B(a,b,.): [0,\infty)\rightarrow [0,1]$ is continuous,
\item[xii.] $lim_{\lambda\rightarrow \infty}B(a,b,\lambda)=0$ \quad $(\forall\lambda >0)$,
\item [xiii.] $Y(a,b,\lambda)=0$ \quad (for $\lambda >0$)  if and only if $a=b$,
\item[xiv.] $Y(a,b,\lambda)=Y(b,a,\lambda)$ \quad $(\forall\lambda >0)$,
\item[xv.] $Y(a,b,\lambda)\bullet Y(b,c,\mu)\geq Y(a,c,\lambda+\mu)$ \quad $(\forall \lambda, \mu >0)$,
\item[xvi.] $Y(a,b,.): [0,\infty)\rightarrow [0,1]$ is continuous,
\item[xvii.] $lim_{\lambda \rightarrow \infty}Y(a,b,\lambda)=0$ \quad (for $\lambda >0$),
\item[xviii.] If $\lambda \leq 0$, then $G(a,b,\lambda)=0$, $B(a,b,\lambda)=1$ and $Y(a,b,\lambda)=1$.
\end{itemize}
Then $\mathcal{N}=(G,B,Y)$ is called Neutrosophic metric(NM) on $F$.
\end{definition}
The functions $G(a,b,\lambda), B(a,b,\lambda), Y(a,b,\lambda)$ denote the degree of nearness, the degree of neutralness and the degree of non-nearness between $a$ and $b$ with respect to $\lambda$, respectively.

\begin{definition}\cite{KirSim}
Give $V$ be a NMS, $0<\varepsilon<1$, $\lambda >0$ and $a\in F$. The set
$O(a,\varepsilon,\lambda)=\lbrace b\in F: G(a,b,\lambda)> 1-\varepsilon, \quad  B(a,b,\lambda)< \varepsilon, \quad Y(a,b,\lambda)< \varepsilon\rbrace$ is said to be the open ball (OB) (center $a$ and radius $\varepsilon$ with respect to $\lambda$).
\end{definition}

\begin{lemma}\cite{KirSim}
Every OB $O(a,\varepsilon,\lambda)$ is an open set (OS).
\end{lemma}


\section{Fixed Point Results}

\begin{definition}
Let $F$ be a set. A non-negative real-valued function $h$ on $F\times F$ is called as a quasi-metric on $F$ if it satisfies the following axioms:
\begin{itemize}
\item [i.] $h(a,b)=h(b,a)=0$ if and only if $a=b$,
\item [ii.] $h(a,b)\leq h(a,c)+h(c,b)$,
\end{itemize}
for all $a,b,c\in F$.
\end{definition}

From this definition we can understand: It is possible $h(a,b)\neq h(b,a)$ for some $a,b\in F$.\\

A quasi-metric is a distance function which satisfies the triangle inequality but is not symmetric in
general. Quasi-metrics are a subject of comprehensive investigation both in pure and applied mathematics
in areas such as in functional analysis, topology and computer science.

\begin{proposition}
Let $V$ be the NMS. For any $\varepsilon \in (0,1]$, define $h: F \times F \rightarrow R^{+}$ as follows:
\begin{eqnarray}\label{fixeq:1}
h_{\varepsilon}(a,b)=\inf\lbrace \lambda>0 : G(a,b, \lambda)> 1-\varepsilon, \quad B(a,b, \lambda)<\varepsilon, \quad Y(a,b, \lambda)<\varepsilon\rbrace
\end{eqnarray}
Then,
\begin{itemize}
\item [i.] $(F, h_{\varepsilon}: \varepsilon \in (0,1])$ is a generating space of quasi-metric family.
\item [ii.] The topology $\tau_{\mathcal{N}}$ on $(F, h_{\varepsilon}: \varepsilon \in (0,1])$ coincides with the $\mathcal{N}-$topology on $V$, that is, $h_{\varepsilon}$ is a compatible symmetric for  $\tau_{\mathcal{N}}$.
\end{itemize}
\end{proposition}

\begin{proof}
Firstly, we prove that (i.). It can be easily seen that $h_{\varepsilon}$ suffices the conditions of the definition of quasi-metric. Let's show the condition (ii.) of quasi-metric. We know that the operations $\circ, \bullet$ are continuous. If we consider Remark \ref{remark:1}, for any given $\varepsilon \in (0,1)$, we can take $\varepsilon^{*} \in (0, \varepsilon)$ such that $(1-\varepsilon^{*}) \circ (1-\varepsilon^{*})> 1-\varepsilon$ and $\varepsilon^{*} \bullet \varepsilon^{*}< \varepsilon$. Given $h_{\varepsilon}(a,b)=x$ and $h_{\varepsilon}(b,c)=y$. From (\ref{fixeq:1}),
\begin{eqnarray*}
G(a,b, x+\lambda)> 1-\varepsilon^{*}, \quad B(a,b, x+\lambda)<\varepsilon^{*}, \quad Y(a,b, x+\lambda)<\varepsilon^{*}
\end{eqnarray*}
and
\begin{eqnarray*}
G(a,c, y+\lambda)> 1-\varepsilon^{*}, \quad B(a,c, y+\lambda)<\varepsilon^{*}, \quad Y(a,c, y+\lambda)<\varepsilon^{*}.
\end{eqnarray*}
From here,
\begin{eqnarray*}
G(a,c, x+y+2\lambda)\geq G(a,b, x+\lambda) \circ G(b,c, y+\lambda) > (1-\varepsilon^{*}) \circ (1-\varepsilon^{*}) > 1-\varepsilon,
\end{eqnarray*}

\begin{eqnarray*}
B(a,c, x+y+2\lambda)\leq  B(a, b, x+\lambda) \bullet  B(b,c, y+\lambda) < \varepsilon^{*} \bullet \varepsilon^{*} < \varepsilon,
\end{eqnarray*}

and

\begin{eqnarray*}
Y(a,c, x+y+2\lambda)\leq  Y(a, b, x+\lambda) \bullet  Y(b,c, y+\lambda) < \varepsilon^{*} \bullet \varepsilon^{*} < \varepsilon.
\end{eqnarray*}

Therefore, we have $h_{\varepsilon}(a,c) \leq x+y+2\lambda = h_{\varepsilon}(a,b)+h_{\varepsilon}(b,c)+2\lambda$. Since $\lambda>0$ is arbitrary, $h_{\varepsilon}(a,c) \leq  h_{\varepsilon}(a,b)+h_{\varepsilon}(b,c)$.\\

Now, we prove that (ii.). We must show that
\begin{eqnarray*}
h_{\varepsilon}(a,c)<\lambda \quad \Leftrightarrow \quad G(a,b, \lambda)> 1-\varepsilon, \quad B(a,b, \lambda)<\varepsilon, \quad Y(a,b, \lambda)<\varepsilon
\end{eqnarray*}
for any $\lambda>0$ and $\varepsilon \in (0,1)$. If $h_{\varepsilon}(a,b)<\lambda$, then $G(a,b, \lambda)> 1-\varepsilon, \quad B(a,b, \lambda)<\varepsilon, \quad Y(a,b, \lambda)<\varepsilon$ from (\ref{fixeq:1}).\\

Conversely, consider $G(a,b, \lambda)> 1-\varepsilon, \quad B(a,b, \lambda)<\varepsilon, \quad Y(a,b, \lambda)<\varepsilon$. Since the functions $G,B,Y$ are continuous from the Definition, then there exists an $\eta>0$ such that $G(a,b, \lambda-\eta)> 1-\varepsilon, \quad B(a,b, \lambda-\eta)<\varepsilon, \quad Y(a,b, \lambda-\eta)<\varepsilon$. From here, we have $h_{\varepsilon}(a,b)\leq \lambda-\eta < \lambda$.
\end{proof}

\begin{definition}
Let $V$ be a NMS. The mapping $f:F\rightarrow F$ is called neutrosophic contraction(NC) if there exists $k \in (0,1)$ such that
\begin{eqnarray*}
\frac{1}{G(f(a), f(b), \lambda)}-1 \leq k(\frac{1}{G(a,b,\lambda)-1}), \quad B(f(a), f(b), \lambda)\leq k B(a,b,\lambda), \quad
Y(f(a), f(b), \lambda)\leq k Y(a,b,\lambda)
\end{eqnarray*}
for each $a,b\in F$ and $\lambda >0$.
\end{definition}

\begin{definition}
Let $V$ be a NMS and let $f:F\rightarrow F$ be a NC mapping. Then there exists $c\in F$ such that $c=f(c)$. That is, $c$ is called neutrosophic fixed point (NFP) of $f$.
\end{definition}

Generally, we claim that the contractions have fixed point. If all contractions(including NC) have fixed points, then we can easily say that
$f^{2}$ should have a fixed point. In below proposition, we will show that if $f^{n}$ is a
NC then, $f^{n}$ has fixed point.\\

\begin{proposition}\label{fppro:1}
Suppose that $f$ is a NC. Then $f^{n}$ is also a NC. Furthermore, if $k$ is the constant for $f$, then $k^{n}$ is the constant for $f^{n}$.
\end{proposition}

\begin{proof}
We will use the induction for proof. We take $n=2$. If $f$ is a NC, then, it is clear that
\begin{eqnarray}\label{fpe:1}
h(f(x), f(y))\leq k \times h(x,y)
\end{eqnarray}
for $k\in (0,1)$. If we apply $f$ to both of sides of inequality (\ref{fpe:1}), we have
\begin{eqnarray}\label{fpe:2}
h(f^{2}(a), f^{2}(b))\leq k \times h(f(a),f(b)).
\end{eqnarray}
From (\ref{fpe:1}), we can write $k \times h(f(a),f(b)) \leq k^{2} \times h(a,b)$. Thus, if we combine the (\ref{fpe:1}) and (\ref{fpe:2}) and last inequality, we have
\begin{eqnarray*}
h(f^{2}(a), f^{2}(b))\leq k \times h(f(a),f(b)) \leq k^{2} \times h(a,b)
\end{eqnarray*}
which leads us to the fact that $f^{2}$ is a NC. If we consider that $f^{n}$ is a NC, then, we can say that $f^{n+1}$ is a NC with above processes.\\

We must prove that the constant for $f^{n}$ is $k^{n}$. Consider $h(f^{n}(a), f^{n}(b))\leq  k^{n} \times h(a,b)$. As similar to the above process, we can apply $f$ to both of sides of this inequality, we get
\begin{eqnarray*}
h(f^{n+1}(a), f^{n+1}(b))\leq k^{n+1} \times h(f(a),f(b)) \leq k^{n+1} \times h(a,b).
\end{eqnarray*}
From this inequality, $h(f^{n+1}(a), f^{n+1}(b))\leq k^{n+1} \times h(a,b)$. Therefore, we understand that the Theorem is true for all $n$.
\end{proof}

\begin{remark}
From Proposition \ref{fppro:1}, we can say that each $f^{n}$ has the same fixed point. Because, if we take $f(a)=a$, then $f^{2}=f(f(a))=f(a)=a$ and by induction, $f^{n}(a)=a$.
\end{remark}

\begin{proposition}\label{fppro:2}
Let $f$ be a NC and $a\in F$. $f[O(a,\varepsilon, \lambda)]\subset O(a,\varepsilon, \lambda)$ for large enough values of $\varepsilon$.
\end{proposition}

\begin{proof}
Let $b\in O(a,\varepsilon, \lambda)$. We must find $\varepsilon$ such that $f(b) \in O(a,\varepsilon, \lambda)$ and so $h(a, f(b))<\varepsilon$. We can write
\begin{eqnarray*}
h(a, f(b))\leq h(a, f(a))+h(f(a), f(b)).
\end{eqnarray*}
Further, $h(f(a), f(b))\leq k \times h(a,b)$ and $h(a,b)\leq \varepsilon$. Therefore, $h(f(a), f(b))\leq k \times \varepsilon$. From this inequalities, we have $h(a, f(b))\leq h(a, f(a))+k \times \varepsilon$. We can choose $\varepsilon$ so that $h(a, f(a))+k \times \varepsilon$. Then, for any $b\in O(a,\varepsilon, \lambda)$, $h(x, f(b))<\varepsilon$ and $f(b)\in O(a,\varepsilon, \lambda)$.
\end{proof}

\begin{remark}
From Proposition \ref{fppro:2} and the definitions neutrosophic open ball and neutrosophic closed ball, if the inclusion $f[O(a, \varepsilon, \lambda)]\subset O(a, \varepsilon, \lambda)$ is hold, then the inclusion also $\overline{f[O(a, \varepsilon, \lambda)]}\subset \overline{O(a, \varepsilon, \lambda)}$ is hold.
\end{remark}

\begin{proposition}\label{fppro:3}
The inclusion $f^{n}[O(a, \varepsilon, \lambda)]\subset O(f^{n}(a), \varepsilon^{*}, \lambda)$ is hold for all $n$, where $\varepsilon^{*}=k^{n}\times \varepsilon$.
\end{proposition}

The proof of this proposition is similar to Proposition \ref{fppro:1}.

\begin{remark}
It is fact that if the inclusion $f^{n}[O(a, \varepsilon, \lambda)]\subset O(f^{n}(a), \varepsilon^{*}, \lambda)$ is hold, then the inclusion also $\overline{f^{n}[O(a, \varepsilon, \lambda)]}\subset \overline{O(f^{n}(a), \varepsilon^{*}, \lambda)}$ is hold.
\end{remark}

\begin{theorem}\label{fpthe:1}
Let $V$ be a complete NMS. Let $f:F \rightarrow F$ be a NC mapping. Then, $f$ has a unique NFP.
\end{theorem}

Theorem \ref{fpthe:1} is a consequence of Theorem 3.6 in \cite{RaNo}. Hence, using the consept of neutrosophy, Theorem \ref{fpthe:1} is proved as similar Theorem 3.6 in \cite{RaNo}.\\

For the alternative proof, we can use $\overline{f^{n}(O(a, \varepsilon, \lambda))}$. That is, if we choose $b\in \overline{f^{n}(O(a, \varepsilon, \lambda))}$, we can see that $f(b)\in \overline{f^{n}(O(a, \varepsilon, \lambda))}$. Therefore, the distance between $b$ and $f(b)$ is $\epsilon$ and so $f(b)=b$. Thus $b$ is a fixed point.

\section{Conclusion}

The purpose of this paper is to apply the NMS which defined by Kirisci and Simsek \cite{KirSim}.
NC mapping is defined. After the properties related to NC are proved, fixed point theorem
is given.


\end{document}